\documentclass[a4paper,12pt]{amsart}
\usepackage{amssymb}
\usepackage{graphicx}
\usepackage{ifthen}
\nonstopmode \numberwithin{equation}{section}
\newtheorem{thm}{Theorem}
\newtheorem{cor}{Corollary}
\newtheorem{lem}{Lemma}

\newtheorem{conj}{Conjecture}

\theoremstyle{definition}
\newtheorem{defn}{Definition}[section]

\newtheorem{prob}[equation]{Problem}

\newenvironment{rem}{%
\bigskip
\noindent \textsl{{\sl Remark. }}}{\bigskip}
\newenvironment{rems}{%
\bigskip
\noindent \textsl{{\sl Remarks. }}}{\bigskip}

\newcounter {own}
\def\theown {\thesection       .\arabic{own}}

\newenvironment{pf}[1][]{%
 \vskip 3mm
 \noindent
 \ifthenelse{\equal{#1}{}}%
  {{\slshape Proof. }}%
  {{\slshape #1.} }%
 }%
{\qed\bigskip}

\newcounter{alphabet}
\newcounter{tmp}
\newenvironment{Thm}[1][]{\refstepcounter{alphabet}%
\bigskip%
\noindent%
{\bf Theorem \Alph{alphabet}}%
\ifthenelse{\equal{#1}{}}{}{ (#1)}%
{\bf .} \itshape}{\vskip 8pt}

\newenvironment{Lem}[1][]{\refstepcounter{alphabet}%
\bigskip%
\noindent%
{\bf Lemma \Alph{alphabet}}%
\ifthenelse{\equal{#1}{}}{}{ (#1)}%
{\bf .} \itshape}{\vskip 8pt}

\newcommand{\ID}{{\mathbb D}}
\newcommand{\IN}{{\mathbb N}}
\newcommand{\IC}{{\mathbb C}}
\newcommand{\K}{{\mathcal K}}

\newcommand{\M}{{\mathcal M}}

\def\be{\begin{equation}}
\def\ee{\end{equation}}

\newcommand{\bee}{\begin{enumerate}}
\newcommand{\eee}{\end{enumerate}}

\newcommand{\blem}{\begin{lem}}
\newcommand{\elem}{\end{lem}}
\newcommand{\bthm}{\begin{thm}}
\newcommand{\ethm}{\end{thm}}
\newcommand{\bcor}{\begin{cor}}
\newcommand{\ecor}{\end{cor}}
\newcommand{\beg}{\begin{examp}}
\newcommand{\eeg}{\end{examp}}
\newcommand{\begs}{\begin{examples}}
\newcommand{\eegs}{\end{examples}}
\newcommand{\bdefe}{\begin{defn}}
\newcommand{\edefe}{\end{defn}}
\newcommand{\bprob}{\begin{prob}}
\newcommand{\eprob}{\end{prob}}
\newcommand{\bei}{\begin{itemize}}
\newcommand{\eei}{\end{itemize}}

\newcommand{\bcon}{\begin{conj}}
\newcommand{\econ}{\end{conj}}
\newcommand{\bcons}{\begin{conjs}}
\newcommand{\econs}{\end{conjs}}
\newcommand{\bprop}{\begin{propo}}
\newcommand{\eprop}{\end{propo}}
\newcommand{\br}{\begin{rem}}
\newcommand{\er}{\end{rem}}
\newcommand{\brs}{\begin{rems}}
\newcommand{\ers}{\end{rems}}
\newcommand{\bo}{\begin{obser}}
\newcommand{\eo}{\end{obser}}
\newcommand{\bos}{\begin{obsers}}
\newcommand{\eos}{\end{obsers}}
\newcommand{\bpf}{\begin{pf}}
\newcommand{\epf}{\end{pf}}
\newcommand{\ba}{\begin{array}}
\newcommand{\ea}{\end{array}}
\newcommand{\beq}{\begin{eqnarray}}
\newcommand{\beqq}{\begin{eqnarray*}}
\newcommand{\eeq}{\end{eqnarray}}
\newcommand{\eeqq}{\end{eqnarray*}}

\newcommand{\ov}{\overline}

\newcounter{minutes}
\divide\time by 60
\newcounter{hours}
\multiply\time by 60 \addtocounter{minutes}{-\time}
\begin{document}
\title{A note on the Bohr inequality}
\begin{center}
{\tiny \texttt{FILE:~\jobname .tex,
        printed: \number\year-\number\month-\number\day,
        \thehours.\ifnum\theminutes<10{0}\fi\theminutes}
}
\end{center}
\author{Bappaditya Bhowmik${}^{\mathbf{*}}$}
\address{Bappaditya Bhowmik, Department of Mathematics,
Indian Institute of Technology Kharagpur, Kharagpur - 721302, India.}
\email{bappaditya@maths.iitkgp.ac.in}
\author{Nilanjan Das}
\address{Nilanjan Das, Department of Mathematics,
Indian Institute of Technology Kharagpur, Kharagpur - 721302, India.}
\email{nilanjan@iitkgp.ac.in}

\subjclass[2010]{30A10, 30B10, 30C45, 30C80, 30H05}
\keywords{Bohr radius, Analytic functions, Subordination.\newline
${}^{\mathbf{*}}$ Corresponding author}

\begin{abstract}
This article focuses on the Bohr radius problem for
the derivatives of analytic functions, along with a technique of establishing Bohr
inequalities in classical and generalized settings.
\end{abstract}
\thanks{The first author of this article would like to thank the Council of Scientific and Industrial Research
(CSIR), India (Ref.No.- 25(0281)/18/EMR-II) for its financial support.}

\maketitle
\pagestyle{myheadings}
\markboth{B. Bhowmik, N. Das}{A note on the Bohr inequality}

\bigskip
\section{Introduction}
Given any complex valued analytic function $F(z)=\sum_{n=0}^\infty A_nz^n$ defined in
the open unit disk $\ID$, let the \textit{majorant series} of $F$ be
$\M_F(r):=\sum_{n=0}^\infty |A_n|r^n$, $|z|=r<1$.
Using this notation, we state the following remarkable result from \cite{Bohr}.

\begin{Thm}\label{P1TheoA}
Let $f(z)= \sum_{n=0}^{\infty}a_{n}z^n$ be analytic in $\ID$ and $|f(z)|<1$ for all $z\in\ID$. Then
$\M_f(r)\leq 1$ for all $z\in\ID$ with $|z|=r\leq1/3$.
This radius $1/3$ is the best possible.
\end{Thm}
This theorem was initially proved by Harald Bohr for $r\leq 1/6$, and was refined thereafter by
Wiener, Riesz and Schur independently.
Problems of similar type are being extensively studied nowadays in different
settings (see for instance \cite{Aiz, Boas, Def, Kohr, Paul, Pop} and the references therein), and
have gained popularity by the name \textit{Bohr phenomenon}. In view of Theorem~A, it is natural to ask if the derivative of an analytic self
map $f$ of $\ID$ also has Bohr phenomenon. Taking $f(z)=f_a(z)=(z+a)/(1+az)$, $z\in\ID$ for each $a\in[0,1)$, we observe that
$f_a^\prime(z)=(1-a^2)/(1+az)^2$, and hence $\M_{f^\prime_a}(r)\leq 1$ if and only if $|z|=r\leq r^\prime=a/(1+\sqrt{1-a^2})$.
Since $r^\prime\to 0$ as $a\to 0$, it is impossible to obtain a Bohr radius independent of $f$ for $f^\prime$
without additional assumption.
In this article a study of the Bohr inequalities for derivatives is initiated. Namely,
we prove that an analogue of Theorem~A exists for $f^\prime$ if $f:\ID\to\ID$ is an analytic function
with $f(0)=0$.

We now introduce the concepts of \textit{subordination} and \textit{majorization} to facilitate our discussion.
For two analytic functions $f$ and $g$ in $\ID$,
we say $g$ is subordinate to $f$ if there exists a function $\phi$, analytic
in $\ID$ with $\phi(0)=0$ and $|\phi(z)|<1$, satisfying $g=f\circ\phi$.
Throughout this article we denote $g$ is subordinate to $f$  by $g\prec f$. Also
the class of functions $g$ subordinate to a fixed function $f$ will be denoted by $S(f)$.
Again, for two analytic functions $f$ and $g$ in $\ID$, $g$ is said to be majorized by $f$ if $|g(z)|\leq|f(z)|$ for all $z\in\ID$.
Besides having wide application in Geometric function theory, these two concepts are
closely connected to the very rich theory of composition and multiplication operators.
A lot of research from function theoretic aspect has been devoted to establish ``comparison theorems'' between $|g^\prime(z)|$ and
$|f^\prime(z)|$, provided that $g$ is subordinate to or majorized by $f$ (see \cite[pp. 206-210]{Dur} and \cite{Mac}). However,
the superordinate function $f$ has always been chosen to be univalent or locally univalent in all such results.
In fact, local univalence is necessary for such results to hold. To see this, let $|g^\prime(z)|\leq|f^\prime(z)|$ for all
$|z|\leq\alpha$, $\alpha>0$ whenever $g\prec f$, $f$ not necessarily locally univalent on all of $\ID$. Now
for any $0<\alpha_1<\min\{1,\alpha\}$, we choose $f(z)=(z-\alpha_1)^2$
and $g(z)=(z^2-\alpha_1)^2\prec f(z), z\in\ID$. Hence $|g^\prime(z)|=|4z(z^2-\alpha_1)|\leq|f^\prime(z)|=|2(z-\alpha_1)|$ for all $|z|\leq\alpha_1.$
Clearly $g^\prime(\alpha_1)=0$, which is only possible if $\alpha_1=0$ or $1$. If the assumption $g\prec f$ is replaced by
$|g(z)|\leq|f(z)|$ for all $z\in\ID$, then we need to choose $g(z)=z^3-2\alpha_1z^2$ and $f(z)=z^2-2\alpha_1z$
(cf. \cite[p. 95]{Mac}), and argue as above.
The second theorem of this paper shows that such comparison results will hold for $\M_{g^\prime}(r)$ and $\M_{f^\prime}(r)$
with $g$ and $f$ being merely analytic.

At this point, it should be mentioned that comparison results for $\M_g(r)$ and $\M_f(r)$ are already available in literature,
provided that $g\prec f$ or $g$ is majorized by $f$. One of them is
the following lemma from \cite{BB}, which is extended to a more general case in \cite[Theorem 2.1]{Alk}.

\begin{Lem}\label{P4lemA}
If $g(z)=\sum_{n=0}^\infty b_nz^n\prec f(z)=\sum_{n=0}^\infty a_nz^n, z\in\ID$ then $\M_g(r)\leq \M_f(r)$
for all $|z|=r\leq 1/3$.
\end{Lem}
This radius $1/3$ in Lemma~B is pointed out to be the best possible in \cite[Corollary 2.2]{Alk}. To prove that,
one can consider $f(z)=z$ and
$$
g(z)=\xi_a(z):=z(z-a)/(1-az)\prec f(z)
$$
for each $a\in[0,1)$. From this we have $\M_g(r)\leq \M_f(r)$ for $|z|=r\leq 1/(1+2a)$,
and hence letting $a\to 1-$ shows that the radius $1/3$ in Lemma~B cannot be improved in general.
Now as an answer to the problem of determining the Bohr radius for the subordinating
family of an odd analytic function, it was shown in \cite[Theorem 2.4]{Alk} that the
conclusion of Lemma~B holds for $|z|=r\leq 1/\sqrt{3}$ if $f$ and $g$ both are taken to be odd analytic. However,
since $f(z)=z$ is odd, the example we just provided shows that the radius $1/3$ cannot be improved for $g\in S(f)$ if only $f$ is odd analytic.
The same example also shows that the radius $1/3$ in \cite[Corollary 2.3]{Alk} cannot be
improved if any analytic function $g$ is majorized by an odd analytic function
$f$. In this paper we have shown that an improvement on $1/3$
is possible in this situation if $f$ and $g$ both are odd analytic.
Further, we prove another comparison theorem between $\M_{g^\prime}(r)$ and $\M_{f^\prime}(r)$, assuming that
$f$ and $g$ both are odd analytic, and that $g\prec f$ or $g$ is majorized by $f$ in $\ID$.

Now we discuss a different approach towards the Bohr inequality.
In the single variable framework, it is a natural problem to
find an appropriate form of the Bohr inequality for the functions which map
$\ID$ inside a domain in $\IC$ other than $\ID$.
To this end, a generalized treatment of the Bohr radius problem has been introduced in \cite{Abu}, using the concept of subordination.
According to \cite{Abu},
we say that $S(f)$ has the Bohr phenomenon if for any $g(z)=\sum_{n=0}^{\infty}b_{n}z^n\in S(f)$, there is a $r_0\in(0,1]$
such that
\be\label{P4eq4}
\M_g(r)-|b_0|\leq d(f(0),\partial f(\ID))
\ee
for $|z|=r<r_0$. Here $d(f(0),\partial f(\ID))$ denotes the Euclidean distance between $f(0)$ and the boundary of domain $f(\ID)$.
It is seen that whenever an analytic function $g$ maps $\ID$ into a domain $\Omega$ other than $\ID$, then in a general sense the
Bohr inequality $(\ref{P4eq4})$ can be established if $g$ can be recognized as a member of $S(f)$, $f$ being the covering map
from $\ID$ onto $\Omega$ satisfying $f(0)=g(0)$.
In particular, if we take $\Omega=\ID$, then for any analytic $g:\ID\to\Omega$ there exists a disk automorphism $f$
such that $g(0)=f(0)$ and $g\in S(f)$.
In this case $ d(f(0),\partial\ID)=1-|f(0)|$, and hence $(\ref{P4eq4})$ reduces to the classical Bohr inequality $\M_g(r)\leq 1$.
The Bohr phenomenon has been researched by using the above definition in a number of articles
(see for example \cite{Abu, Abu1, Ali, BB, BB1} and their references).
In view of this approach, it is often necessary to
obtain an estimate on $\M_g(r)$ in terms of $\M_f(r)$ for establishing the Bohr inequality $(\ref{P4eq4})$,
and the aforesaid Lemma~B is already seen to be useful for this purpose (cf. \cite{BB}).
Clearly, any comparison result of this type with a different radius
might come in handy in several situations.
We prove another lemma in this regard and discuss its interesting outcomes, which include
a refinement of \cite[Theorem 2.7]{Alk}, and the Bohr phenomena for \textit{spherically convex functions}
and \textit{convex functions of bounded type}. These two concepts will be explained in due course.
\section{Main results}
We start by stating two important properties of the majorant series, which will be used
frequently in a number of our proofs. Let $f(z)=\sum_{n=0}^\infty a_nz^n$ and $g(z)=\sum_{n=0}^\infty b_nz^n$ be
two analytic functions defined in $\ID$. Then
$\M_{f+g}(r)\leq \M_f(r)+\M_g(r)$ and $\M_{fg}(r)\leq \M_f(r)\M_g(r)$ for any $|z|=r\in[0,1)$.
For the sake of completeness, we include the proofs for both of these properties.
Indeed, it is easy to see that $\M_{f+g}(r)=\sum_{n=0}^\infty |a_n+b_n|r^n\leq \sum_{n=0}^\infty |a_n|r^n+\sum_{n=0}^\infty |b_n|r^n
=\M_f(r)+\M_g(r)$. Further, using the same argument we can prove that if $f_c(z)=\sum_{i\in\mathbb{Z}}f_i(z)$ is analytic in
$\ID$ where $f_i$ is analytic in $\ID$ for each $i\in\mathbb{Z}$, then $\M_{f_c}(r)\leq \sum_{i\in\mathbb{Z}}\M_{f_i}(r)$.
To prove the other property, we only need to note the trivial facts $\M_{\alpha f}(r)=|\alpha|\M_f(r)$ for any $\alpha\in\IC$
and $\M_{z^t f}(r)= r^t\M_f(r)$. Now $(fg)(z)=\sum_{n=0}^\infty a_n(z^n g(z))$, which readily gives
$\M_{fg}(r)\leq\sum_{n=0}^\infty |a_n|r^n\M_g(r)=\M_f(r)\M_g(r).$

We are now ready to
establish the Bohr phenomenon for the derivative of an analytic self map $w$ of $\ID$,
where $w(0)=0$.
\bthm\label{P4thm6}
Let $w(z)=\sum_{n=1}^\infty w_nz^n$ be an analytic self map of $\ID$. Then
$\M_{w^\prime}(r)\leq 1$ for $|z|=r\leq r_0=1-\sqrt{2/3}$. This radius $r_0$ is the best possible.
\ethm
\bpf
By the Schwarz lemma, $w(z)=z\phi(z)$, where $\phi:\ID\to\ov{\ID}$ is an analytic function. Therefore
$w^\prime(z)=z\phi^\prime(z)+\phi(z)$, which yields $\M_{w^\prime}(r)\leq r\M_{\phi^\prime}(r)+\M_\phi(r)=
|w_1|+\sum_{n=2}^\infty n|w_n|r^{n-1}$. Using
Wiener's estimates $|w_n|\leq 1-|w_1|^2$, $n\geq 2$, we have $\M_{w^\prime}(r)\leq |w_1|+2(1-|w_1|)\sum_{n=2}^\infty nr^{n-1}$.
As a result, a little computation reveals that $\M_{w^\prime}(r)\leq 1$ if $2(1/(1-r)^2-1)\leq 1$, i.e. if $r\leq r_0=1-\sqrt{2/3}$.
To see that $r_0$ is the best possible
we consider the functions
$w(z)=\xi_a(z), z\in\ID$, where $a\in[0,1)$.
For any fixed $a\in[0,1)$,
$$
\M_{\xi^\prime_a}(r)=a+(1-a^2)\left(\frac{r}{1-ar}+\frac{r}{(1-ar)^2}\right)\leq 1
$$
if and only if
$$
M(a,r):=(1+a)\left(\frac{r}{1-ar}+\frac{r}{(1-ar)^2}\right)-1\leq 0,
$$
which holds if and only if $r\leq r_0(a)$, $r_0(a)$ being a nonnegative real number. 
Now we observe that for each fixed $a$, $M(a,r)$ is strictly increasing in
$r\in (0,1)$, and that $M(a,0)=-1<0$, $M(a,1)=(1+3a-2a^2)/(1-a)^2>0$.
Hence the number $r_0(a)$ is the only real root of $M(a,r)=0$ in $(0,1)$. Again, we observe that $M(a,r)$ is strictly increasing in $a\in[0,1)$
for each fixed $r\in(0,1)$, and therefore  for any $a_1, a_2\in[0,1)$,
$r_0(a_1)<r_0(a_2)$ whenever $a_1>a_2$. As a consequence, $M(a,r)\leq 0$ for all $a\in[0,1)$
if and only if $r\leq\inf_{a\in[0,1)}r_0(a)$, and $\inf_{a\in[0,1)}r_0(a)=\lim_{a\to 1-}r_0(a)$ is the root of $\lim_{a\to 1-}M(a,r)=0$.
Noting that $\lim_{a\to 1-}M(a,r)=2(1/(1-r)^2-1)-1$, we conclude $\inf_{a\in[0,1)}r_0(a)=1-\sqrt{2/3}=r_0$. This completes the proof.
\epf

The above theorem allows us to derive the following result.
\bthm\label{P4thm7}
Let $f(z)=\sum_{n=0}^\infty a_nz^n$ and $g(z)=\sum_{n=0}^\infty b_nz^n$ be two analytic functions defined in $\ID$. Then for
$|z|=r\leq r_0=1-\sqrt{2/3}$
\bee
\item[(i)]~
if $g\prec f$ then $\M_{g^\prime}(r)\leq \M_{f^\prime}(r)$
\item[(ii)]~
if $|g(z)|\leq|f(z)|$ for all $z\in\ID$, and $f(0)=0$ then $\M_{g^\prime}(r)\leq \M_{f^\prime}(r)$.
\eee
The radius $r_0$ cannot be improved in either case.
\ethm
\bpf
(i) By definition, $g\prec f$ implies $g(z)=f(w(z))$, $z\in\ID$ where $w$ is an analytic self map of $\ID$ with $w(0)=0$.
Since $g^\prime(z)=w^\prime(z)f^\prime(w(z))$, $\M_{g^\prime}(r)\leq \M_{w^\prime}(r)\M_{f^\prime\circ w}(r)$. As $f^\prime\circ w\prec f^\prime$,
from Lemma~B we have $\M_{f^\prime\circ w}(r)\leq \M_{f^\prime}(r)$ for $r\leq 1/3$; and from Theorem \ref{P4thm6}, $\M_{w^\prime}(r)\leq 1$
for $r\leq r_0=1-\sqrt{2/3}<1/3$.
Hence $\M_{g^\prime}(r)\leq \M_{f^\prime}(r)$ for $r\leq r_0=1-\sqrt{2/3}$.

\noindent
(ii) If $|g(z)|=|f(z)|$ for some $z\in\ID,$ then by the maximum modulus principle, $g(z)=cf(z)$ for some $|c|=1$, which implies
$\M_{g^\prime}(r)=\M_{f^\prime}(r)$ for all $r\in[0,1)$. Therefore assuming $|g(z)|<|f(z)|$ for all $z\in\ID$, there exists an analytic self map
$\phi$ of $\ID$ such that $g(z)=\phi(z)f(z)$, $z\in\ID$. Now $g^\prime(z)=\phi^\prime(z)f(z)+\phi(z)f^\prime(z)$. Further, observing that
$\M_{f/z}(r)\leq \M_{f^\prime}(r)$ for all $r\in[0,1)$, we have
$\M_{g^\prime}(r)\leq \M_{z\phi^\prime}(r)\M_{f/z}(r)+\M_\phi(r)\M_{f^\prime}(r)\leq (\M_{z\phi^\prime}(r)+\M_\phi(r))\M_{f^\prime}(r)$.
Following similar lines of calculations as in the proof of Theorem \ref{P4thm6} it could be shown that
$\M_{z\phi^\prime}(r)+\M_\phi(r)=r\M_{\phi^\prime}(r)+\M_\phi(r)\leq 1$ for $r\leq 1-\sqrt{2/3}$,
and hence $\M_{g^\prime}(r)\leq \M_{f^\prime}(r)$ for
$r\leq r_0= 1-\sqrt{2/3}$.

To see that this $r_0$ is the best possible, one can choose the functions
$g(z)=\xi_a(z), z\in\ID$, where $a\in[0,1)$ and $f(z)=z$. Clearly $\xi_a\prec f$ and $|\xi_a(z)|\leq |f(z)|$ for all $z\in\ID$.
Since $f^\prime(z)=1$, the argument used to establish the best possible part in the proof of Theorem \ref{P4thm6} will also
show that the radius $r_0$ cannot be improved in both (i) and (ii).
\epf

\br
The condition $f(0)=0$ is essential in Theorem \ref{P4thm7}(ii). For, if we consider $f(z)=z+1$ and $g(z)=z(z+1)$, $z\in\ID$
then $|g(z)|\leq |f(z)|$ in $\ID$, but $\M_{g^\prime}(r)=2r+1>1=\M_{f^\prime}(r)$ for all $r\neq 0$.
\er

Our next result exhibits that $\M_g(r)$ is dominated by $\M_f(r)$ for $|z|=r\leq 1/\sqrt{3}$, provided that $f$ and $g$ both
are odd analytic and $g$ is majorized by $f$ in $\ID$.
\bthm\label{P4thm1}
Let $f(z)=\sum_{n=0}^\infty a_{2n+1}z^{2n+1}$ and $g(z)=\sum_{n=0}^\infty b_{2n+1}z^{2n+1}$ be two odd analytic functions defined in $\ID$ such that
$|g(z)|\leq |f(z)|$ for all $z\in\ID$. Then $\M_g(r)\leq \M_f(r)$ for all $|z|=r\leq 1/\sqrt{3}$. This radius $1/\sqrt{3}$ is the best possible.
\ethm
\bpf
If $|g(z)|=|f(z)|$ for any $z\in\ID$, then by the maximum modulus principle $g(z)=cf(z)$ for some $|c|=1$ and for all
$z\in\ID$. This implies $\M_g(r)=\M_f(r)$ for all $|z|=r<1$. Therefore considering $|g(z)|<|f(z)|$ for all $z\in\ID$, we see that
there exists an analytic $\phi:\ID\to\ID$ with a Taylor expansion $\phi(z)=\sum_{n=0}^\infty c_nz^n$, $z\in\ID$ such that $g(z)=\phi(z)f(z)$ for all $z\in\ID$.
As the coefficients for $z^{2n}$ for any $n\geq 0$ is zero for both $f$ and $g$, it is immediately seen that
$$
\sum_{n=0}^\infty b_{2n+1}z^{2n+1}=\left(\sum_{n=0}^\infty c_{2n}z^{2n}\right)\left(\sum_{n=0}^\infty a_{2n+1}z^{2n+1}\right), z\in\ID.
$$
Hence for $|z|=r\in[0,1)$, $\M_g(r)\leq \M_{\hat{\phi}}(r)\M_f(r)$, where $\hat{\phi}(z):=\sum_{n=0}^\infty c_{2n}z^{2n}$.
Using $|c_{2n}|\leq 1-|c_0|^2$ for all $n\geq 1$,  we get that $\M_{\hat{\phi}}(r)\leq 1$ for
$r\leq 1/\sqrt{3}$, and hence $\M_g(r)\leq \M_f(r)$ for $r\leq 1/\sqrt{3}$. To see that $1/\sqrt{3}$ is the best possible,
we take
$$
g(z)=g_a(z):=\xi_a(z^2)/z, z\in\ID\,
$$
for each $a\in[0,1)$ and $f(z)=z$. Clearly $|g_a(z)|\leq|f(z)|$ for all $z\in\ID$ and for any fixed $a\in[0,1)$,
and $\M_{g_a}(r)\leq\M_f(r)$ for $|z|=r\leq 1/\sqrt{1+2a}$.
Now allowing $a\to 1-$, our proof will be complete.
\epf

\br
Using the same functions $g_a$ and $f$ from the proof of Theorem \ref{P4thm1} one can show that the radius $1/\sqrt{3}$ is the best
possible in \cite[Theorem 2.4]{Alk} as well.
\er

The following theorem reveals that the radius $1-\sqrt{2/3}$ in Theorem \ref{P4thm7} can also be improved if $f$ and $g$ both are taken
to be odd analytic.
\bthm\label{P4thm9}
Let $f(z)=\sum_{n=0}^\infty a_{2n+1}z^{2n+1}$ and $g(z)=\sum_{n=0}^\infty b_{2n+1}z^{2n+1}$ be two odd analytic functions defined in $\ID$. Then for
$|z|=r\leq \tilde{r_0}=\sqrt{(4-\sqrt{13})/3}$
\bee
\item[(i)]~
if $g\prec f$ then $\M_{g^\prime}(r)\leq \M_{f^\prime}(r)$
\item[(ii)]~
if $|g(z)|\leq|f(z)|$ for all $z\in\ID$, then $\M_{g^\prime}(r)\leq \M_{f^\prime}(r)$.
\eee
This radius $\tilde{r_0}$ cannot be improved in either case.
\ethm
\bpf
(i) We have $g(z)=f(w(z))$ where $w$ is an analytic self map of $\ID$ with $w(0)=0$.
As $f, g$ both are odd analytic, comparing the coefficients from both sides of
the equation $\sum_{n=0}^\infty b_{2n+1}z^{2n+1}=\sum_{n=0}^\infty a_{2n+1}{w}^{2n+1}(z)$ gives
$g(z)=\sum_{n=0}^\infty a_{2n+1}z^{2n+1}\phi_n(z)$, where $\phi_n(z)=(w^{2n+1}(z)-w^{2n+1}(-z))/2z^{2n+1}$
is an even analytic function satisfying $\phi_n(\ID)\subset{\ov{\ID}}$. Now
$g^\prime(z)=\sum_{n=0}^\infty a_{2n+1}(z^{2n+1}\phi_n^\prime(z)+(2n+1)z^{2n}\phi_n(z))$, and hence
$\M_{g^\prime}(r)\leq\sum_{n=0}^\infty |a_{2n+1}|(r^{2n+1}\M_{\phi_n^\prime}(r)+(2n+1)r^{2n}\M_{\phi_n}(r))
\leq\sum_{n=0}^\infty (2n+1)|a_{2n+1}|r^{2n}(r\M_{\phi_n^\prime}(r)+\M_{\phi_n}(r)).$
Let $\phi_n(z)=\sum_{k=0}^\infty \phi_n^kz^{2k}.$ Using the Wiener's estimates $|\phi_n^k|\leq 1-|\phi_n^0|^2$ for each $k\geq1$,
a little calculation yields that
$r\M_{\phi_n^\prime}(r)+\M_{\phi_n}(r)=\sum_{k=0}^\infty (2k+1)|\phi_n^k|r^{2k}\leq 1$ whenever
$2\sum_{k=1}^\infty (2k+1)r^{2k}\leq1 $, i.e. whenever $r\leq\tilde{r_0}$. Thus $\M_{g^\prime}(r)\leq
\sum_{n=0}^\infty(2n+1)|a_{2n+1}|r^{2n}=\M_{f^\prime}(r)$ for $r\leq\tilde{r_0}$.

\noindent
(ii) This part can be proved by adopting exactly the same lines of computations from the proof of Theorem \ref{P4thm7}(ii).
We only need to note that $\phi$ can be assumed to be even analytic (see the proof of Theorem \ref{P4thm1}), and therefore
proceeding like the proof of part (i), we can show that
$\M_{z\phi^\prime}(r)+\M_\phi(r)=r\M_{\phi^\prime}(r)+\M_\phi(r)\leq 1$ for $r\leq\tilde{r_0}$.

To see that $\tilde{r_0}$ is the best possible in both (i) and (ii), one can choose $g(z)=g_a(z), z\in\ID$ where $a\in[0,1)$
and $f(z)=z$. It has to be shown that $\M_{g_a^\prime}(r)\leq\M_{f^\prime}(r)=1$ for all $a\in[0,1)$ if and only if
$r\leq\tilde{r_0}$. The argument hereafter is mostly similar to the corresponding part from
the proof of Theorem \ref{P4thm6}.
\epf

We now prove a lemma which will be instrumental in proving the forthcoming results.
\blem\label{P4lem1}
If $g(z)=\sum_{n=m}^\infty b_nz^n\prec f(z)=\sum_{n=m}^\infty a_nz^n, z\in\ID$ with $a_n\neq 0$ for all $n\geq m$, $m\in\IN$, then $\M_g(r)\leq \M_f(r)$
for all $|z|=r\leq \inf_{n\geq m}|a_{n+1}/a_n|$. If $m=0$,
$\M_g(r)\leq\M_f(r)$ for all $|z|=r\leq \inf_{n\geq 1}|a_{n+1}/a_n|$.
\elem
\bpf
From the proof of \cite[Theorem 6.3]{Dur} it is easy to see that, whenever $g\prec f$ and $m\in\IN$,
$$
\sum_{n=m}^t \frac{|b_n|^2}{|a_n|}r^n\leq \sum_{n=m}^t|a_n|r^n
$$
if $\{r^n/|a_n|\}_{n=m}^t$ is a monotonic decreasing sequence, i.e. if $r\leq |a_{n+1}/a_n|$ for all $n$ satisfying $m\leq n\leq t$.
Therefore for $r\leq  \inf_{m\leq n\leq t}|a_{n+1}/a_n|$, using the Cauchy-Schwarz inequality gives
$$
\sum_{n=m}^t|b_n|r^n\leq \sqrt{\left(\sum_{n=m}^t \frac{|b_n|^2}{|a_n|}r^n\right)}
\sqrt{\left(\sum_{n=m}^t|a_n|r^n\right)}\leq \sum_{n=m}^t |a_n|r^n.
$$
Letting $t\to\infty$, the proof for the first part is complete. For $m=0$,
the result follows from the fact that $g(z)-b_0\prec f(z)-a_0$ for all $z\in\ID$.
\epf

\br
It is evident that for $f(z)=\sum_{n=m}^\infty a_nz^n, z\in\ID$ and $m\in\IN$,
if $\{|a_n|\}_{n=m}^t$ is a monotonic increasing sequence with $a_m\neq 0$, then
for any $g(z)=\sum_{n=m}^\infty b_nz^n\prec f(z)$ we have $\sum_{n=m}^t |b_n|r^n\leq \sum_{n=m}^t |a_n|r^n$ for all $|z|=r<1$. This
can be taken as a majorant series analogue of the much esteemed result \cite[Theorem VII(ii)]{Rogo}.
\er

The following result recovers \cite[Theorem 2.7]{Alk}, with the constant $2\sqrt{3}-3\approx 0.4641016$ replaced
by a much smaller number.
\bthm\label{P4thm8}
Let $f:\ID\to\ID$ be an analytic function with Taylor expansion $f(z)=\sum_{n=0}^\infty a_nz^n$. Then
the inequality $|f(z)|+\M_f(r)-|a_0|\leq 1$ holds for all $|z|=r\leq r_{|a_0|}$, where
$$
r_{|a_0|}:=\frac{\sqrt{(1+|a_0|)^2+|a_0|^2}-(1+|a_0|)}{|a_0|^2}
$$
and $|a_0|\geq\tilde{a}\approx 0.361103$, $\tilde{a}$ being the only
root in $(0,1)$ of the equation $x^4+2x^2+2x-1=0$. This radius $r_{|a_0|}$
is sharp.
\ethm
\bpf
We assume $a_0\neq 0$ to begin with. Now
using the fact that $f(z)\prec (z+a_0)/(1+\ov{a_0}z), z\in\ID$, from Lemma \ref{P4lem1}
we get that
$$
\M_f(r)-|a_0|\leq \frac{r(1-|a_0|^2)}{1-|a_0|r}
$$
for $|z|=r\leq |a_0|$. Note that this result is
originally due to Bombieri (see \cite[Theorem 8.6.1]{Ga}). As our method differs from Bombieri's own,
it could be considered as an alternative proof.
Now to establish the theorem, one
has to follow exactly the same lines of the proof of \cite[Theorem 2.7]{Alk}, only the computations are
to be carried out under the restriction $r\leq |a_0|$ instead of $r\leq 1/3$. Hence
$|f(z)|+\M_f(r)-|a_0|\leq 1$ for $r\leq r_{|a_0|}$ if there exists $\tilde{a}\in(0,1)$
such that $r_{|a_0|}\leq |a_0|\iff |a_0|\geq \tilde{a}.$
Set $|a_0|=x.$ Looking at $\psi(x)=(x^3+x+1)-\sqrt{(1+x)^2+x^2}$, we observe that $\psi(x)\geq 0$ if and only if
$(x^3+x+1)^2\geq(1+x)^2+x^2$, which is equivalent to $\psi_1(x)=x^4+2x^2+2x-1\geq 0$. Now since $\psi_1(0)=-1<0$
and $\psi_1(1)=4>0$, $\psi_1$ has a root in $(0,1)$ which
we choose to be $\tilde{a}$. Since $\psi_1^\prime(x)>0$ for all $x\in(0,1)$, $\psi_1(x)\geq 0\iff x\geq \tilde{a}$. Also,
$\psi_1(2\sqrt{3}-3)=908-524\sqrt{3}>0$, which implies $\tilde{a}<2\sqrt{3}-3$. Using Mathematica, one can see that
$\tilde{a}\approx 0.361103$.
\epf

Before we proceed further we need to be familiar with the following two classes of functions.
A meromorphic univalent (i.e. one-one) function $f$ defined in $\ID$ is called \textit{spherically convex} if $f(\ID)$ is a spherically
convex domain in the Riemann sphere $\hat{\IC}$, i.e. for any $w_1, w_2\in f(\ID),$ the
smaller arc of the greatest circle (spherical geodesic) between
$w_1$ and $w_2$ lies in $f(\ID)$.
Let $\K_s(\alpha)$ be the class of all spherically convex functions
in the normalized form $f(z)=\alpha z+a_2z^2+a_3z^3+\cdots$, $0<\alpha\leq 1.$
The reader is referred to the articles \cite{Mej, Min} and the references therein for a
detailed study on spherically convex functions.
Now for an analytic locally univalent function $f$ in $\ID$ the curvature of $f(\{z:|z|=r\})$ at the point
$f(z)$ is given by
$$
\kappa_f(z)= \frac{1}{|zf^\prime(z)|}\mbox{Re}\left(1+z\frac{f^{\prime\prime}(z)}{f^\prime(z)}\right), 0<|z|=r<1.
$$
A convex analytic univalent function $f$ defined in $\ID$ is said to be
in $CV(R_1, R_2)$ if $f$ satisfies $f(0)=f^\prime(0)-1=0$ and
$$
0<R_1\leq \underline{\lim}_{|z|\to 1-}(1/\kappa_f(z))\leq \overline{\lim}_{|z|\to 1-}(1/\kappa_f(z))\leq R_2<\infty.
$$
The class $CV(R_1, R_2)$ was introduced in \cite{Good} and functions in this class are called \textit{convex functions of
bounded type}. More on $CV(R_1, R_2)$ can be found in \cite{Good1, Wirths}. To the best of our knowledge, coefficient problem
in any of these two classes is not completely solved yet, which makes the Bohr radius problem worth trying. Applying Lemma \ref{P4lem1} enables
us to do the same in the next two theorems.
For brevity we denote $d(f(0), \partial f(\ID))$ by $\delta$ from now onwards.
\bthm\label{P4thm4}
Let $f(z)=\alpha z+\sum_{n=2}^\infty a_nz^n\in\K_s(\alpha)$ with $0<\alpha\leq\sqrt{3}/2$. Then
$\M_f(r)\leq\delta$ for all $|z|=r\leq r_s=1/(1+2\sqrt{1-\alpha^2})$. This radius $r_s$ is sharp for
the function $k_\alpha(z)=\alpha z/(1-\sqrt{1-\alpha^2}z), z\in\ID$.
\ethm
\bpf
From \cite[Corollary 2, p. 132]{Min}, it follows that $\delta\geq \alpha/(1+\sqrt{1-\alpha^2})$.
Now from the inequality $(2.6)$ of \cite[Theorem 3]{Mej}, we observe that $\K_s(1)=\{z\}$, and for $\alpha<1$,
$f(z)=\alpha z/(1+ z\phi(z)\sqrt{1-\alpha^2}), z\in\ID$ where $\phi:\ID\to\ov{\ID}$ is an analytic function.
This is equivalent to saying that $f(z)/\alpha z\prec 1/(1+z\sqrt{1-\alpha^2})$, and therefore a use of Lemma \ref{P4lem1}
gives $(1/\alpha)(\sum_{n=2}^\infty|a_n|r^{n-1})\leq r\sqrt{1-\alpha^2}/(1-r\sqrt{1-\alpha^2})$ for
$|z|=r\leq \sqrt{1-\alpha^2}$. Simplifying a little bit, we get that for $|z|=r\leq \sqrt{1-\alpha^2}$
$$
\M_f(r)=\alpha r+\sum_{n=2}^\infty |a_n|r^n\leq \frac{\alpha r}{1-r\sqrt{1-\alpha^2}}.
$$
Now to have $\M_f(r)\leq\delta$ satisfied, we need
$\alpha r/(1-r\sqrt{1-\alpha^2})\leq\alpha/(1+\sqrt{1-\alpha^2})$, which happens whenever $r\leq r_s= 1/(1+2\sqrt{1-\alpha^2})$.
Note that $r_s\leq\sqrt{1-\alpha^2}$ if and only if $\alpha\leq \sqrt{3}/2$, and hence the proof is done.
The sharpness part can be verified from direct computations.
\epf

\bthm\label{P4thm5}
Let $f(z)=z+\sum_{n=2}^\infty a_nz^n\in CV(R_1, R_2)$ with $R_2\geq 2\delta$. Then
$\M_f(r)\leq\delta$ for all $|z|=r\leq r_c=R_2/(3R_2-2\delta)$. This radius is sharp
for the function $k_a(z)=z/(1-az), z\in\ID$ where $a\in [1/2,1)$.
\ethm
\bpf
We have $f(z)\prec Bz/(1-Az)$ where $A= (R_2-\delta)/R_2$ and $B=(2R_2-\delta)\delta/R_2$ (cf. \cite[p. 544]{Good}).
Under the assumption $R_2\geq 2\delta$,
both $A$ and $B$ are strictly positive. Therefore a use of Lemma \ref{P4lem1} yields
$r+\sum_{n=2}^\infty |a_n|r^n\leq Br/(1-Ar)$ for $|z|=r\leq A$. A little computation reveals that $Br/(1-Ar)\leq \delta$ whenever
$r\leq \delta/(B+A\delta)=r_c$, and that $r_c\leq A$ for $R_2\geq 2\delta$. Now we observe that for any $k_a(z), a\in(0,1)$ defined in $\ID$,
$\delta=1/(1+a)$, $R_2=1/(1-a^2)$, and therefore the sharpness of this result
can be proved from straightforward calculations. Also, $R_2\geq 2\delta$ for $k_a$ if and only if $a\geq 1/2$.
\epf

\brs
\bee
\item
From \cite[Theorem 2]{Good} we have $\delta\geq R_2-\sqrt{R_2^2-R_2}$,
and therefore $r_c\geq R_2/(R_2+2\sqrt{R_2^2-R_2})$.
Hence it is possible to obtain a Bohr radius $r_c$ in the above theorem which is independent of $\delta$.
\item
One can replace the hypothesis $R_2\geq 2\delta$ in the statement of Theorem \ref{P4thm5} by $\delta\in (1/2, 2/3]$, as from the Corollary in \cite[p. 544]{Good}
it is immediate that for $1/2<\delta\leq 2/3$, $R_2\geq \delta^2/(2\delta-1)\geq 2\delta$.
\eee
\ers
It is worth noting that both $r_s$ and $r_c$ turn out to be bigger or equal to $1/3$.
As a result, Lemma~B is not of much use in either case, while using Lemma
\ref{P4lem1} solves the Bohr radius problem to certain extent.

\end{document}